\documentclass[12pt]{article}

\usepackage{t1enc,amsmath,amssymb,amsthm,amscd,epsfig,latexsym,amsfonts,ifthen}

\usepackage{amsmath}%
\usepackage{amsfonts}%
\usepackage{amssymb}%
\usepackage{graphicx}
\newtheorem{theorem}{Theorem}
\theoremstyle{plain}

\newtheorem{definition}{Definition}

\newtheorem{lemma}{Lemma}

\newtheorem{proposition}{Proposition}
\newtheorem{remark}{Remark}

\numberwithin{equation}{section}
\begin{document}

\title{\vspace{-1in}\parbox{\linewidth}{\footnotesize\noindent}
 \vspace{\bigskipamount} \\
Operators associated to the Cauchy-Riemann operator in elliptic complex numbers
\thanks{ 
\hfil\break\indent
{\em $\dagger$ danieldaniel@gmail.com ~
\em $\ddagger$ cvanegas@usb.ve }}}

\author{Daniel Alay\'on-Solarz$\dagger$ and Carmen Judith Vanegas$\ddagger$ \\
Departamento de Matem\'aticas Puras y Aplicadas \\
Universidad Sim\'{o}n Bol\'{\i}var, Caracas 1080-A, Venezuela}

\date{}

\maketitle

\begin{abstract}

In this article we provide a generalized version of the result of L.H. Son and W. Tutschke \cite{tut} on the solvability of first order systems on the plane whose initial functions are arbitrary holomorphic functions. This is achieved by considering the more general concept of holomorphicity with respect to the structure polynomial $X^2 + \beta X + \alpha$. It is shown that the Son-Tutschke lemma on the construction of complex linear operators associated to the Cauchy-Riemann operator remains valid when interpreted for a large class of real parameters $\alpha$ and $\beta$ including the elliptic case but also cases that are not elliptic. For the elliptic case, first order interior estimates are obtained via the generalized version of the Cauchy representation theorem for elliptic numbers and thus the method of associated operators is applied to solve initial value problems with initial functions that are holomorphic in elliptic complex numbers. 

\end{abstract}

\maketitle

\section{Introduction} 
We consider linear first order systems 
\begin{equation}\label{unouno}
\partial_{t} u = a_{11} \partial_{x} u + a_{12} \partial_{y} u + a_{21} \partial_{x} v + a_{22} \partial_{y} v + c_{1}  u + c_{2}  v + c_{3}  
\end{equation}
\begin{equation}\label{unodos}
\partial_{t} v = b_{11} \partial_{x} u + b_{12} \partial_{y} u + b_{21} \partial_{x} v + b_{22} \partial_{y} v + d_{1}  u + d_{2}  v + d_{3}  
\end{equation}
for two real-valued functions $u(t,x,y)$ and $v(t,x,y)$ where $t$ is the time and $(x,y)$ runs in a bounded domain in the $x,y$-plane. The coefficients are supposed to depend at least continuosly on $t,x$ and $y$. Let us recall that in view of the Cauchy-Kovaleskaya Theorem the initial value problem
\begin{equation}\label{unotres}
u(0,x,y) = \phi(x,y)
\end{equation}
\begin{equation}\label{unocuatro}
v(0,x,y) = \psi(x,y)
\end{equation}
is solvable provided the coefficients of (\ref{unouno}) and (\ref{unodos}) and the initial functions posess power-series representation. In our main reference \cite{tut} there are formulated conditions on the coefficients of (\ref{unouno}) and (\ref{unodos}) under which each initial value problem (\ref{unotres}) and (\ref{unocuatro}) is solvable with the additional condition that the initial functions $\phi$ and $\psi$ satisfy the Cauchy-Riemann conditions. Moreover it is shown that under certain conditions four coefficients of (\ref{unouno}) and (\ref{unodos}) may not posess power-series representation, in fact, remarkably they can be chosen as arbitrary continuous functions. 

The goal of this article is to generalize the study the solvability of the system (\ref{unouno}) and (\ref{unodos}) when the initial value functions are not necessarily holomorphic in the ordinary sense but instead they are holomorphic in the more general sense of algebras with structure polynomial, as in the language developed in \cite{Tut2}. For the plane, we consider the structure polynomial $X^2 + \beta X + \alpha$, where $\alpha$ and $\beta$ are real numbers. In these algebras a complex number is of the form $z = x + i y$ where
\begin{displaymath}
i^2 = -\beta i - \alpha
\end{displaymath}
and two real function $u(x,y)$ and $v(x,y)$ satisfy the Cauchy-Riemann equations if
\begin{displaymath}
 \partial_x u - \alpha \partial_y v = 0 ~~ \mbox{and} ~~ \partial_y u + \partial_x v - \beta \partial_y v = 0.
\end{displaymath}

A key tool for the main result in \cite{tut} and also here is the following; we briefly recall that 
a pair of differentiable operators \textbf{L} and \textbf{G} is said to be $\textit{associated}$ if $\textbf{G}(w) = 0$ 
implies $\textbf{G} (\textbf{L} (w)) = 0$, that is if $\textbf{L}$ sends null-solutions of $\textbf{G}$ again 
into null-solutions of $\textbf{G}$. Let
\begin{displaymath}
\textbf{L}w : = A \partial_z w + B \overline{\partial_z w} + C \partial_{\bar z} w + D \overline{\partial_{\bar z} w}
+ Ew + F\overline{w} + G.
\end{displaymath}
where $A,B,E,F$ and $G$ are differentiable functions on $z$ and $\overline z$. Then:
\begin{proposition}( Son-Tutschke lemma)
Suppose  $A,B,E,F$ and $G$ are continuosly differentiable with respect to $z$ and $\overline z$. 
Then the operator \textbf{L} is associated to the Cauchy-Riemann operator if and only if the following conditions are satisfied:
\begin{itemize}
\item B and F are identically equal to zero.
\item A, E and G are holomorphic.
\end{itemize}
\end{proposition}
We will see that this result remains valid for a large family of parameters $\alpha$ and $\beta$ if 
we interpret products as well as the Cauchy-Riemann and conjugated Cauchy-Riemann operators with respect 
to the structure polynomial $X^2 + \beta X + \alpha$. For the elliptic case, that is when $4 \alpha - \beta^2 > 0$ 
in view of the Cauchy representation formula first order interior estimates can be obtained for the derivative of a 
holomorphic function and thus an abstract Cauchy-Kovaleskaya theorem is applicable.

This article is structured as follows; in the section \textbf{2} we introduce the complex numbers with structure polynomials $X^2 + \beta X + \alpha$ as well as the concept of holomorphic functions in these algebras and we recall the Cauchy-Pompeiu and Cauchy representation formulas that hold for the elliptic case. In the section \textbf{3} we introduce the complex derivative for the elliptic case and provide some elementary properties. In the section \textbf{4} we find first order interior estimates for the complex derivative in the elliptic case. In the section \textbf{5} we find different operators  for a large class of structure polynomials for which the Cauchy-Riemann operator is associated. 
In the section \textbf{6} we provide the complex rewriting of the system (\ref{unouno}) and (\ref{unodos}) 
for the structure polynomial $X^2 + \beta X + \alpha$. In the section \textbf{7} we characterize the real system (\ref{unouno}) and (\ref{unodos}). 
In the section \textbf{8} we prove the Weirstrass approximation theorem for elliptic numbers and finally, in the section \textbf{9} we give our main result on the solvability of the initial value problem of the said system.

\section{Preliminaries} 

Let us consider the complex numbers $z = x + iy$, endowed with the parameter-depending complex algebra with structure polynomial $X^2 + \beta X + \alpha$, that is where $i^2 = -\alpha - \beta i$, and $\alpha$ and $\beta$ are real numbers. For two numbers $z_1 = x_1 + i y_1$ and $z_2 = x_2 + i y_2$ the product induced is given by
\begin{displaymath}
z_1 z_2 = (x_1 x_2 - \alpha  y_1 y_2) + i (x_1 y_2 + x_2 y_1 - \beta  y_1 y_2).
\end{displaymath}
It is easy to show that for all real $\alpha$ and $\beta$ this product is both commutative and associative. When $\alpha = 1$ and $\beta =0$ this reduces to the ordinary complex numbers, if  $\alpha = -1$ and $\beta =0$ this algebra reduces to the ordinary hyperbolic numbers (or also \textit{double numbers}) and if $\alpha = \beta = 0$ this algebra reduces to the ordinary parabolic numbers (or also \textit{dual numbers}). The conjugation, denoted $\bar z$ is given by  $x - iy$.  

The Cauchy-Riemann operator $\partial_{\bar{z}}$ is defined by
$ \partial_{\bar{z}} = \frac{1}{2} (\partial_{x} + i \partial_{y}) $
and a function $w(z) = u(x,y) + i v(x,y)$ is said to be holomorphic in this algebra  if $\partial_{\bar{z}} w = 0$.  Holomorphicity of $w = u + iv$ with respect to $X^2 + \beta X + \alpha$  is equivalent to $u$ and $v$ satisfying the following system of equations:
\begin{displaymath}
 \partial_x u - \alpha \partial_y v = 0, ~~~
\partial_y u + \partial_x v - \beta \partial_y v = 0,
\end{displaymath}
which is a generalization of the ordinary Cauchy-Riemann equations. The conjugated Cauchy-Riemann operator $\partial_z$ is given by
$ \partial_{z} = \frac{1}{2} (\partial_{x} - i \partial_{y}). $
The Cauchy-Riemann operator and its conjugate commute with each other, so for any function $w$ we have
$ \partial_{z} \partial_{\bar{z}} w = \partial_{\bar{z}} \partial_{z} w.$
The second order operator $\partial_{z} \partial_{\bar{z}}$ applied to a function $w = u + i v$ and equalizing 
to zero yields the following system of equations:
\begin{displaymath}
\frac{1}{4}(\partial_{xx} u + \alpha  \partial_{yy} u - \alpha \beta   \partial_{yy} v) = 0, ~~~ 
\frac{1}{4}(\beta  \partial_{yy} u + \partial_{xx} v + (\alpha - \beta^2 )  \partial_{yy} v  ) = 0,
\end{displaymath}
which is uncoupled if and only if  $\beta = 0$.

A complex number $z =  x + iy$ admits an inverse if $ x^2 - \beta  xy + \alpha  y^2 \neq 0$, in this case the reciprocal map is given by
\begin{displaymath}
(x + iy)^{-1} := \frac{x - \beta y - i y}{ x^2 - \beta  xy + \alpha  y^2 }.
\end{displaymath}
If  $\alpha$ and $\beta$ satisfy the ellipticity condition $4 \alpha - \beta^2 > 0$ then every number different from zero is invertible as the polynomial $x^2 - \beta xy + \alpha  y^2$ vanishes only at the origin. For the elliptic case the square root of this  polynomial defines a norm depending on $\alpha$ and $\beta$:
\begin{displaymath}
|z|_{(\alpha, \beta)} = \sqrt{x^2 - \beta x y + \alpha y^2}.
\end{displaymath}
This norm is compatible with the algebra with structure polynomial $X^2 + \beta X + \alpha$  in the sense that for all $z_1$ and $z_2$
\begin{displaymath}
|z_1 \cdot z_2 |_{(\alpha, \beta)} = |z_1  |_{(\alpha, \beta)} | z_2 |_{(\alpha, \beta)}
\end{displaymath}
and for all $z \neq 0$
\begin{displaymath}
|z^{-1}|_{(\alpha, \beta)} = \frac{1}{|z|_{(\alpha, \beta)}}.
\end{displaymath}
The euclidian norm is $|\cdot |_{(1,0)}$ which we denote by $| \cdot |$. As all norms over finite-dimensional spaces are equivalent, for $\alpha$ and $\beta$ satisfying the ellipticity condition numbers
$K_1(\alpha, \beta)$ and $K_2(\alpha, \beta)$ exist such that for all $z$ we have
\begin{displaymath}
K_1(\alpha, \beta) |z|_{(\alpha, \beta)} \leq |z| \leq K_2(\alpha, \beta) |z|_{(\alpha, \beta)}.
\end{displaymath}

Many formulas from the ordinary complex analysis remain the same in the generalized case, for example the product rule
for two differentiable functions $f, g$:
\begin{displaymath}
\partial_{ \overline z} (fg)= f \cdot \partial_{ \overline z}  g + \partial_{ \overline z} f \cdot g.
\end{displaymath}
Other results are only valid for the elliptic case, and are written slightly different, for example the Cauchy-Pompeiu formula as proved in \cite{AV} is given by
\begin{displaymath}
f(\zeta) = \frac{1}{ 2 \pi \hat{i} }  \int\limits_{\partial \Omega} \frac{f(z)}{\widetilde{z - \zeta}} d \tilde z - \frac{1}{ \pi \hat{i}}  \iint\limits_{ \Omega} \frac{\partial_{\bar z} f(z)}{\widetilde{z - \zeta}} dx dy 
\end{displaymath}
where $f$ is a continuously differentiable function on $ \bar{\Omega} $,~ $\Omega$ is a simple connected and bounded domain
with sufficiently smooth boundary,
$\tilde{z} = y - ix $,  $d \tilde{z} = dy - i dx $ and $ \hat{i} = \frac{\beta + 2 i}{\sqrt{4 \alpha - \beta^2 }} $.
In particular if the function $f$ is holomorphic, for the elliptic case one gets the Cauchy integral representation formula
\begin{displaymath}
f(\zeta) = \frac{1}{ 2 \pi \hat{i} }  \int\limits_{\partial \Omega} \frac{f(z)}{\widetilde{z - \zeta}} d \tilde z. 
\end{displaymath}

\section{The complex derivative}

As illustrated in the introduction, certain formulae are to be written slightly different than in the ordinary case. One example being the Cauchy-Pompeiu representation formula. Another example is the notion of the complex derivative. It is convenient to define the complex derivative for the elliptic case in the following manner:
\begin{definition}
Let $f : \Omega \subseteq \mathbb{C}(\alpha, \beta) \to \mathbb{C}(\alpha, \beta)$ be a complex function. We say $f$ is \textit{differentiable} at $z_0 \in \Omega$ if the following
limit exists:
\begin{displaymath}
\lim_{h \to 0} (-i) \frac{f(z_0+h) - f(z_0)}{\tilde{h}},
\end{displaymath}
where $h = h_1 + i h_2$ and $\tilde{h} = h_2 - i h_1$. If this limit exists, the complex derivative at $z_0$ is denoted as $f'(z_0)$.
\end{definition}
The limit can be taken indistinctly in the euclidian norm or the norm depending on $\alpha$ and $\beta$ as they are topologically equivalent. Since in the ordinary case $\tilde h = i^{-1} h$ we notice that the previous definition coincides with the usual for $\alpha = 1$ and $\beta = 0$. Now suppose $f$ is differentiable at $z_0$, then taking the limit in the direction $h \equiv x$  or $h \equiv i y$ we obtain
\begin{displaymath}
f'(z_0) =  \partial_{x} f (z_0) ~~ \mbox{or} ~~
f'(z_0) = -i \partial_{y} f (z_0),
\end{displaymath}
respectively.
So if the function $f$ is differentiable at $z_0$ then it is a solution to the Cauchy-Riemann equation at this $z_0$:
\begin{displaymath}
\partial_{\bar z} f(z_0) = \frac{1}{2} \big( \partial_{x} f (z_0) + i  \partial_{y} f (z_0) \big) = 0.
\end{displaymath}
On the other side given that
$ \partial_{x} = \partial_{\bar z} + \partial_{z} $
if $f$ is holomorphic (in $\Omega$) then
\begin{equation}\label{fprima}
f'(z) = \partial_{x} f =  \big( \partial_{\bar z} + \partial_{z}  \big) f =  \partial_{z} f.
\end{equation}
So the complex derivative corresponds to the conjugate Cauchy-Riemann operator for an holomorphic function $f$ like in the ordinary case. 
Also, in view of the triangle inequality and for $f = u + iv$:
\begin{displaymath}
|\partial_{x} u | \leq | f'|, ~~~ 
|\partial_{x} v | \leq | f'|.
\end{displaymath}

For the partial derivative depending on $y$ we have:
\begin{displaymath}
 \partial_{y} f  = \frac{\beta + i}{\alpha} f'
\end{displaymath}
then
\begin{displaymath}
\frac{1}{K_2(\alpha, \beta)} | \partial_{y} f  | \leq   | \partial_{y} f  |_{(\alpha, \beta)} = \frac{1}{\sqrt{\alpha}} | f'|_{(\alpha, \beta)} \leq \frac{1}{K_1(\alpha, \beta)\sqrt{\alpha}} | f'|
\end{displaymath}
and thus
\begin{displaymath}
| \partial_{y} f  | \leq \frac{K_2(\alpha, \beta)}{K_1(\alpha, \beta)\sqrt{\alpha}} | f'|.
\end{displaymath}
Then again, in view of the triangle inequality
\begin{displaymath}
|\partial_{y} u | \leq  \frac{K_2(\alpha, \beta)}{K_1(\alpha, \beta)\sqrt{\alpha}}| f'| ~~ \mbox{and} ~~
|\partial_{y} v | \leq \frac{K_2(\alpha, \beta)}{K_1(\alpha, \beta)\sqrt{\alpha}} | f'|.
\end{displaymath}
Like in the ordinary case, the euclidean norm of the complex derivative bounds the euclidean norm of the partial derivatives of the components of the function.

\section{First order interior estimates}

Roughly speaking interior estimates describe the behaviour of the derivatives of a function near the boundary of a bounded domain.
In order to obtain interior estimate for holomorphic functions $f$ we will use 
the Cauchy Integral Formula for $f'$. Let $f$ be holomorphic in the bounded domain $\Omega$ and continuous in 
its closure $\bar{\Omega}$. Consider an arbitrary point $\zeta \in \Omega$. If $r$ is less than the distance 
dist$(\zeta, \partial \Omega)$ of $\zeta$ from the boundary $\partial \Omega$ and in view of the identity
\begin{displaymath}
\partial_{ \xi} \big(\frac{1}{ \widetilde{z - \zeta}} \big) =  - \frac{i}{(\widetilde{z - \zeta})^2}
\end{displaymath}
then applying (\ref{fprima}) the Cauchy Integral Formula for $f'$ implies
\begin{displaymath}
f'(\zeta) = - \frac{i}{ 2 \pi \hat{i} }  \int\limits_{|z - \zeta| = r} \frac{f(z)}{(\widetilde{z - \zeta})^2} d \tilde z.
\end{displaymath}
Taking the norm $| \cdot |_{(\alpha, \beta)}$ we obtain:
\begin{displaymath}
| f'(\zeta) |_{(\alpha, \beta)} = \big|\frac{i}{ 2 \pi \hat{i}}\int\limits_{|z - \zeta| = r} 
\frac{f(z)}{(\widetilde{z - \zeta})^2} d \tilde z \big|_{(\alpha, \beta)} 
\leq \frac{\sqrt{\alpha}}{2 \pi r^2}\sup_{|z - \zeta| = r} | f(z)|_{(\alpha, \beta)} 
\big|\int\limits_{|z - \zeta| = r} d \tilde z\big|_{(\alpha, \beta)}.
\end{displaymath}
Now let us estimate the last factor
$ |\int\limits_{|z - \zeta| = r}  d\tilde z  \ |_{(\alpha, \beta)}$.
Parameterizing $|z - \zeta| = r$ as $z = r( \cos(\theta) + i \sin(\theta) ) + \zeta$ where $0 \leq \theta \leq 2 \pi$ and 
derivating with respect to $\theta$ we obtain
\begin{displaymath}
d z = r( - \sin(\theta) + i \cos(\theta)) \ d \theta ~~~ \mbox{and} ~~~
d \tilde{z} = r ( \cos(\theta) + i \sin(\theta) ) \ d \theta.
\end{displaymath}
Hence
\begin{displaymath}
|  \int\limits_{|z - \zeta| = r}  d\tilde z  \ |_{(\alpha, \beta)} \leq 2 \pi r \sup_{0 \leq \theta \leq \ 2 \pi} | \cos(\theta) + i \sin(\theta) |_{(\alpha, \beta)} \leq \frac{2 \pi r}{K_1(\alpha, \beta)}. 
\end{displaymath}
We then have
\begin{displaymath}
| f'(\zeta) |_{(\alpha, \beta)} \leq \frac{\sqrt{\alpha}}{K_1(\alpha, \beta) r}  \big( \sup_{|z - \zeta| = r} | f(z)|_{(\alpha, \beta)} \big) \leq  \frac{\sqrt{\alpha}}{K^2_1(\alpha, \beta) r}  \sup_{|z - \zeta| = r} | f(z)| 
\end{displaymath}
and given that
\begin{displaymath}
|f' (\zeta) | \leq K_2(\alpha, \beta) | f'(\zeta)|_{(\alpha, \beta)}
\end{displaymath}
we finally obtain
\begin{displaymath}
|f' (\zeta) |  \leq  \frac{K_2(\alpha, \beta)\sqrt{\alpha}}{K^2_1(\alpha, \beta) r}  \big( \sup_{|z - \zeta| = r} | f(z)| \big).
\end{displaymath}
Therefore the limiting process $r \to \mbox{dist}(\zeta, \partial \Omega)$ leads to the interior estimate 
for the derivative of a holomorphic function in terms of the values of the holomorphic function
\begin{displaymath}
|f' (\zeta) |\leq 
\frac{K_2(\alpha, \beta)\sqrt{\alpha}}{K^2_1(\alpha, \beta) \mbox{dist}(\zeta, \partial \Omega)}\big( \sup_{\Omega}| f(z)|\big).
\end{displaymath}
In the ordinary case then $\alpha =1$ and $\beta = 0$ and $K_1(1,0) = K_2(1,0) = 1$, so the latter formula reduces to the well known first order estimate
\begin{displaymath}
|f' (\zeta) |  \leq  \frac{1}{\mbox{dist}(\zeta, \partial \Omega)}  \sup_{\Omega} |f(z)|.
\end{displaymath}

\section{First order differential operators associated to the generalized Cauchy-Riemann operator}

We will prove the following generalization of the Son-Tutschke lemma.
\begin{lemma}
Suppose $\alpha$ and $\beta$ are constants satisfying $\alpha  \beta^2  - 4 \alpha^2 \neq 0$ and $A,B,E,F$ and $G$ are 
continuosly differentiable functions with respect to $z$ and $\overline z$. Then the operator \textbf{L}:
\begin{equation}\label{opeL}
\textbf{L}w : = A \partial_z w + B \overline{\partial_z w} + C \partial_{\bar z} w + D \overline{\partial_{\bar z} w}
+ Ew + F\overline{w} + G.
\end{equation}
is associated to the Cauchy-Riemann operator in the complex algebra with structure polynomial $X^2 + \beta  X + \alpha$ if and only if the following conditions are satisfied:
\begin{itemize}
\item B and F are identically equal to zero.
\item A, E and G are holomorphic.
\end{itemize}
\end{lemma}
At this point we remark that the condition $\alpha  \beta^2  - 4 \alpha^2 \neq 0$ is true for the elliptic case, but also for cases that are not elliptic.
\begin{proof}
Our proof is a suitable modification of the original one. Applying the Cauchy-Riemann operator 
to the operator \textbf{L} defined by (\ref{opeL}), as a consequence of the product rule and 
assuming $w$ is holomorphic with respect to $X + \beta X + \alpha$, we have:
\begin{displaymath}
\partial_{\bar z}A \cdot \partial_z w + \partial_{\bar z}B \cdot \overline{\partial_z w} + B \cdot \partial_{\bar z} (\overline{\partial_z w})+\partial_{\bar z}E \cdot w + \partial_{\bar z}F \cdot \overline{w} + F \cdot \partial_{\bar z}(\overline w) + \partial_{\bar z}G.
\end{displaymath}
So if the conditions on $A, B, E, F$ and $G$ are met then $\textbf{L}w$ is holomorphic 
if $w$ is holomorphic for all $\alpha$ and $\beta$, i.e., $\partial_{\bar{z}}w = 0$ implies $\partial_{\bar{z}}\textbf{L}w = 0$ and
therefore \textbf{L} is associated to $\partial_{\bar{z}}$.

Now assume $\textbf{L}w$ is holomorphic if only $w$ is so and $\alpha  \beta^2  - 4 \alpha^2 \neq 0$. 
In order to obtain the conditions on the coefficients of operator \textbf{L}, 
we will start by choosing special functions of the associated space, in this case
special holomorphic functions, and write the relations saying that \textbf{L} is holomorphic for those functions.
Then choosing  the holomorphic function $w = 0$ we have
$ \partial_{\bar z}G = 0 $
and so $G$ is holomorphic and the term in $G$ can be omitted. We now choose the holomorphic function $w = i$ to get
\begin{displaymath}
\partial_{\bar z}E \cdot i + \partial_{\bar z}F \cdot \overline{i} = (\partial_{\bar z}E  - \partial_{\bar z}F ) \cdot i =0.
\end{displaymath} 
We now observe that $i$ is invertible in the modified complex algebra and
$ i  (-\beta-i) = \alpha, \  \alpha \neq 0 $,
so taking associativity into account we conclude that
\begin{displaymath}
\partial_{\bar z}E  - \partial_{\bar z}F =0.
\end{displaymath}
We now choose the holomorphic function $w = 1$ and we obtain:
\begin{displaymath}
\partial_{\bar z}E  + \partial_{\bar z}F =0
\end{displaymath}
so $E$ and $F$ are holomorphic necessarily.
Now we choose the function $w(x,y)$ $= \alpha x + \beta y + iy$. For this function we get 
$\partial_{\overline{z}} w = 0$, $\partial_z w = \alpha $ and 
$ \partial_{\bar z}(\overline w)= \alpha + i\beta $ to obtain:
\begin{equation}\label{cintro}
 \alpha \cdot \partial_{\bar{z}}A + \alpha \cdot \partial_{\bar{z}}B + F \cdot ( \alpha + \beta  i) =0.
\end{equation}
Now we choose the function $w(x,y)=-y + ix$ which is holomorphic and $\partial_{z} w = i$,  $\partial_{\bar{z}} \overline{ w} = -i$ 
and for this choice we obtain:
\begin{displaymath}
\partial_{\bar{z}}A \cdot i - \partial_{\bar{z}}B \cdot i - F \cdot i = (\partial_{\bar{z}}A - \partial_{\bar{z}}B - F) \cdot i =0.
\end{displaymath}
Since $i$ is invertible we conclude:
\begin{equation}\label{cinco}
\partial_{\bar{z}}A - \partial_{\bar{z}}B - F = 0.
\end{equation}
We now choose the function $w(x,y)= \alpha  x^2 - y^2 +i(\beta  x^2 + 2xy)$, for this function 
we get $\partial_{\overline{z}} w = 0$, $\partial_z w = 2 \alpha  x +2 i (\beta x+y)$, 
$ \partial_{\bar z}(\overline w) =  2 \alpha x - 2i y$ and $\partial_{\bar{z}}( \overline{ \partial_z w}) = 2 \alpha$. So we obtain:
\begin{displaymath}
\partial_{\bar z}A \cdot (2 \alpha x +2 i (\beta  x+y)) + \partial_{\bar z}B \cdot (2 \alpha x - 2 i (\beta x+y)) + 2 \alpha B  + F \cdot ( 2 \alpha x - 2i y) = 0
\end{displaymath}
and after reordering we get the following expression:
\begin{displaymath}
2x \alpha  (\partial_{\bar z}A + \partial_{\bar z}B + F)  + 2i \beta  x (\partial_{\bar z}A -  \partial_{\bar z}B) +2iy ( \partial_{\bar z}A - \partial_{\bar z}B - F ) +2 \alpha B = 0.
\end{displaymath}
Now, taking the equalities (\ref{cintro}) and (\ref{cinco}) into account in the above equality we obtain:
\begin{displaymath}
 2x (- i \beta  F)  + 2ix \beta  F  + 2 \alpha B = 2 \alpha B  = 0,
\end{displaymath}
and thus $B$ vanishes and equations (\ref{cintro}) and (\ref{cinco}) can be rewritten as:
\begin{displaymath}
\alpha   ( \partial_{\bar{z}}A  + F )  + i \beta F  =0 ~~ \mbox{and} ~~ 
\partial_{\bar{z}}A  - F = 0
\end{displaymath}
respectively. 
This in turn yields
$ F  (2 \alpha + i \beta) =0 $.
Since $(2 \alpha + i \beta)(\beta^2 - 2 \alpha + i \beta) =\alpha \beta^2  -4 \alpha^2  \neq 0$, then $2 \alpha + i \beta$ is invertible and
$F$ vanishes and $A$ is holomorphic. 
\end{proof}

\section{The complex rewriting of the given system}

Coming back to the first order system (\ref{unouno}) and (\ref{unodos}) :
\begin{displaymath}
\partial_{t} u = a_{11} \partial_{x} u + a_{12} \partial_{y} u + a_{21} \partial_{x} v + a_{22} \partial_{y} v + c_{1}  u + c_{2}  v + c_{3}  
\end{displaymath}
\begin{displaymath}
\partial_{t} v = b_{11} \partial_{x} u + b_{12} \partial_{y} u + b_{21} \partial_{x} v + b_{22} \partial_{y} v + d_{1}  u + d_{2}  v + d_{3}  
\end{displaymath}
Setting $w = u + i v$, then
\begin{displaymath}
\partial_{\bar{z}}w = \frac{1}{2} (\partial_{x} u - \alpha  \partial_{y} v) + \frac{i}{2}(\partial_{x} v +  \partial_{y} u - \beta \partial_{y} v)
\end{displaymath}
\begin{displaymath}
\overline{\partial_{\bar{z}}w} = \frac{1}{2} (\partial_{x} u - \alpha  \partial_{y} v) - \frac{i}{2}(\partial_{x} v +  \partial_{y} u - \beta \partial_{y} v)
\end{displaymath}
\begin{displaymath}
\partial_{z}w = \frac{1}{2} (\partial_{x} u + \alpha  \partial_{y} v) + \frac{i}{2} (\partial_{x} v -  \partial_{y} u + \beta \partial_{y} v)
\end{displaymath}
\begin{displaymath}
\overline{\partial_{z}w} = \frac{1}{2} (\partial_{x} u + \alpha  \partial_{y} v) - \frac{i}{2} (\partial_{x} v -  \partial_{y} u + \beta \partial_{y} v)
\end{displaymath}

\begin{displaymath}
2 \partial_{x} u = \partial_{z} w + \overline{\partial_{z} w} + \partial_{\bar{z}} w +
\overline{\partial_{\bar{z}} w}
\end{displaymath}
\begin{displaymath}
2 \alpha \partial_{y} v = \partial_{z} w + \overline{\partial_{z} w} - \partial_{\bar{z}} w -
\overline{\partial_{\bar{z}} w}
\end{displaymath}
\begin{displaymath}
2 i  \partial_{x} v = \partial_{z} w - \overline{\partial_{z} w} + \partial_{\bar{z}} w -
\overline{\partial_{\bar{z}} w}
\end{displaymath}
\begin{displaymath}
2 i \alpha \partial_{y} u = (-\alpha + i \beta)\partial_{z} w + (\alpha + i \beta) \overline{\partial_{z} w} + (\alpha -i  \beta) \partial_{\bar{z}} w - (\alpha+ i \beta)
\overline{\partial_{\bar{z}} w}
\end{displaymath}
\begin{displaymath}
2u = w + \overline{w}
\end{displaymath}
\begin{displaymath}
2iv = w - \overline{w}
\end{displaymath}
The complex rewriting of the system (\ref{unouno}) and (\ref{unodos}) is given by
\begin{displaymath}
\partial_t w  = A \partial_z w + B \overline{\partial_z w} + C \partial_{\bar z} w + D \overline{\partial_{\bar z} w}
+ Ew + F\overline{w} + G.
\end{displaymath}
where
\begin{equation}\label{dosA}
2 A = a_{11} + 2 \frac{\beta }{\alpha}  a_{12} - b_{12}  - \frac{\beta}{\alpha} a_{21} +  b_{21} +  \frac{1}{\alpha} a_{22}   + i \big( b_{11} +  \frac{1}{\alpha} a_{12}  + \frac{\beta}{\alpha} b_{12}  - \frac{1}{\alpha} a_{21}   + \frac{1}{\alpha} b_{22} \big)
\end{equation}
\begin{equation}\label{dosB}
2 B = a_{11} + b_{12}  +   \frac{\beta}{\alpha} a_{21}   - b_{21} +  \frac{1}{ \alpha} a_{22} + i \big( b_{11} - \frac{1}{\alpha} a_{12} +  \frac{1}{\alpha} a_{21}  + \frac{\beta}{\alpha} b_{12}  + \frac{1}{\alpha} b_{22}\big)
\end{equation}
\begin{equation}\label{dosC}
2 C = a_{11}  - 2 \frac{\beta}{\alpha} a_{12} + b_{12} - \frac{\beta}{\alpha} a_{21} + b_{21} - \frac{1}{\alpha} a_{22} + i \big( b_{11} - \frac{1}{\alpha} a_{12} - \frac{1}{\alpha} a_{21} - \frac{\beta}{\alpha} b_{12} - \frac{1}{\alpha} b_{22} \big)
\end{equation}
\begin{equation}\label{dosD}
2 D = a_{11} - b_{12} + \frac{\beta}{\alpha} a_{21} - b_{21} - \frac{1}{\alpha} a_{22} + i \big( b_{11} + \frac{1}{\alpha} a_{12} + \frac{1}{\alpha} a_{21} - \frac{\beta}{\alpha} b_{12} - \frac{1}{\alpha} b_{22} \big)
\end{equation}
\begin{equation}\label{dosE}
2 E = c_1 - \frac{\beta}{\alpha} c_2 + d_2 + i \big( -\frac{1}{\alpha} c_2 + d_1 \big)
\end{equation}
\begin{equation}\label{dosF}
2 F = c_1 + \frac{\beta}{\alpha} c_2 - d_2 + i \big( \frac{1}{\alpha} c_2 + d_1 \big)
\end{equation}
\begin{equation}\label{soloG}
G = c_3 + i d_3 
\end{equation}

\section{The characterization of the given real system}

Adding (\ref{dosE}) and (\ref{dosF}) and taking into account that $F = 0$ according to Lemma $1 $, we have
\begin{displaymath}
E = c_1 + i d_1.
\end{displaymath}
From $F = 0$ we obtain
\begin{displaymath}
0 = c_1 + \frac{\beta}{\alpha} c_2 - d_2 + i \big( \frac{1}{\alpha} c_2 + d_1 \big)
\end{displaymath}
and then it follows that
\begin{displaymath}
c_2 = -\alpha d_1, \ \ d_2 = c_1 - \beta d_1
\end{displaymath}
that is, $c_1$, $c_2$, $d_1$ and $d_2$ are determined by the holomorphic function $E$. 
On the other side $c_3$ and $d_3$ are determined by the real and imaginary part of the holomorphic function $G$. 
F
Separating the real and the imaginary parts of the equations (\ref{dosA}) to (\ref{dosD}) we obtain 8 linear equations for the 
eight coefficients $a_{ij}, b_{ij}$, $ i,j=1,2 $: 
\begin{displaymath}
\begin{pmatrix}
1 & 2\frac{\beta}{\alpha} & -1 & -\frac{\beta}{\alpha} & 1 & \frac{1}{\alpha} & 0 & 0 \\
\\
0 & \frac{1}{\alpha} &  \frac{\beta}{\alpha} & -\frac{1}{\alpha} & 0 & 0 & 1 & \frac{1}{\alpha}  \\
\\
1& 0 & 1 & \frac{\beta}{\alpha} & -1 &  \frac{1}{\alpha} & 0 & 0  
\\
\\
0 & - \frac{1}{\alpha} & \frac{\beta}{\alpha} & \frac{1}{\alpha}  & 0 & 0 & 1 & \frac{1}{\alpha}
\\
\\
1 & -2 \frac{\beta}{\alpha} & 1 &  -\frac{\beta}{\alpha} & 1 & -\frac{1}{\alpha} & 0 & 0
\\
\\ 
0 & - \frac{1}{\alpha} & - \frac{\beta}{\alpha} & - \frac{1}{\alpha} & 0 & 0 & 1 & - \frac{1}{\alpha}
\\
\\
1 & 0 & -1 & \frac{\beta}{\alpha} & -1 & - \frac{1}{\alpha} & 0 & 0
\\
\\
0 & \frac{1}{\alpha} & - \frac{\beta}{\alpha} & \frac{1}{\alpha} & 0 & 0 & 1 & -\frac{1}{\alpha}
\end{pmatrix}
\end{displaymath}

The determinant of the matrix of coefficients is $- \frac{256}{\alpha^4} \neq 0$. Therefore, the coefficients $a_{ij}$ and  $b_{ij}$ are uniquely determined if the coefficients $A,B,C$ and $D$ are given. As the coefficients $C$ and $D$ are arbitrary for all the systems associated to the Cauchy-Riemann operator, then 4 of these coefficients can be arbitrarily chosen. Adding (\ref{dosA}) and (\ref{dosB}) 
and recalling that $B=0$ due to Lemma 1, we obtain
\begin{displaymath}
A = A_1 + i A_2 = a_{11} + \frac{\beta}{\alpha} a_{12} + \frac{1}{\alpha} a_{22} + i(b_{11} + \frac{\beta}{\alpha} b_{12} + \frac{1}{\alpha} b_{22})
\end{displaymath}
thus
\begin{displaymath}
a_{22} = \alpha A_1 - \alpha a_{11} - \beta a_{12}
\end{displaymath}
\begin{displaymath}
b_{22} = \alpha A_2 - \alpha b_{11} - \beta b_{12}
\end{displaymath}
where $A = A_1 + i A_2$ is holomorphic according to the Lemma 1. On the other side
subtracting (\ref{dosA}) from (\ref{dosB}) we have
\begin{displaymath}
- A = -A_1 - i A_2 = -\frac{\beta}{\alpha} a_{12} + b_{12} + \frac{\beta}{\alpha} a_{21} - b_{21} - i(\frac{1}{\alpha} a_{12} - \frac{1}{\alpha} a_{21}).
\end{displaymath}
Therefore 
\begin{displaymath}
b_{21} =  A_{1} - \frac{\beta}{\alpha} a_{12} + b_{12} + \frac{\beta}{\alpha} a_{21},
\end{displaymath}
\begin{displaymath}
a_{21} = - \alpha A_2 + a_{12},
\end{displaymath}
combining these two equations yields
\begin{displaymath}
b_{21} = A_1 - \beta A_2 + b_{12}
\end{displaymath}
and therefore, 10 out of the 14 coefficients depend on the choice of 3 arbitrary holomorphic functions (E, G and A), 
the other 4 coefficients are free and it is only imposed the condition that they are continuous functions. We have proved the following Lemma
\begin{lemma}
Suppose the coefficients $a_{11}, a_{12}, b_{11}, b_{12}$ are arbitrarily chosen. Then the coefficients that characterize the system 
(\ref{unouno}) and (\ref{unodos}) are given by:
$c_i$ and $d_i$, $i=1,2,3$ are determined by the choice of 2 arbitrary holomorphic functions ($E$ and $G$). Then permissible coefficients are given by
\begin{displaymath}
a_{21} = -\alpha A_2 + a_{12}
\end{displaymath}
\begin{displaymath}
a_{22} = \alpha A_1 - \alpha a_{11} - \beta a_{12} 
\end{displaymath}
\begin{displaymath}
b_{21} =  A_1 - \beta A_2 + b_{12} 
\end{displaymath}
\begin{displaymath}
b_{22} = \alpha A_2 - \alpha b_{11} - \beta b_{12} \\ \\
\end{displaymath}
where $A_1 + i A_2$ is an arbitrary holomorphic function.
\end{lemma}

\section{The Weirstrass approximation theorem}

Let us consider the path-integral in elliptic complex numbers:
\begin{displaymath}
\int\limits_{C} f(z) \cdot d \tilde{z},
\end{displaymath}
where  $f(z)$ is a continuous function in $C$, a sufficiently smooth path and $d \tilde z = dy - i dx$.  By means of a partition of the smooth path and using the triangle inequality for $| \cdot |_{\alpha, \beta}$ we obtain the estimation:
\begin{displaymath}
| \int\limits_{C} f(z) \cdot d \tilde{z} |_{(\alpha,\beta)} \leq   \int\limits_{C} | f(z) |_{(\alpha,\beta)} \cdot  | d \tilde{z} |_{(\alpha,\beta)},
\end{displaymath}
where $| d \tilde z |_{(\alpha, \beta)}$ is given by
$| d \tilde z |_{(\alpha, \beta)} = \sqrt{\alpha d x^2 + \beta dx dy + d y^2 }$
and thus
\begin{displaymath}
\int\limits_{C} | d \tilde z |_{(\alpha, \beta)}
\end{displaymath}
represents the arc-length of the path considered with respect to the norm depending on $\alpha$ and $\beta$. 

Let $\{ f_n \}_{n = 0}^{\infty}$ be a sequence of continous functions converging to a function $f$ in $\Omega \in \mathbb{C}$, let us suppose, additionally that this convergence is uniform on every compact subset of $\Omega$. With these conditions the function $f$ es continuous in $\Omega$ and for every path $C$ of finite length in $\Omega$ the following estimation holds
\begin{displaymath}
| \int\limits_{C} f(z) \cdot d \tilde{z} - \int\limits_{C} f_{n}(z) \cdot d \tilde{z} |_{(\alpha, \beta)} 
\leq \frac{1}{K^2_1(\alpha, \beta)}\sup_{z \in C} |f(z) - f_n(z)| l(C),
\end{displaymath}
where $l(C)$ denotes the ordinary arc-lenght of $C$. As $C$ is compact then the convergence is uniform. Then for an arbitrary large $n$ the expression in the right hand-side goes to zero. We conclude that
\begin{displaymath}
\lim_{n \to \infty} \int\limits_{C} f_{n}(z) \cdot d \tilde{z} = \int\limits_{C} f(z) \cdot d \tilde{z}. 
\end{displaymath}
Let us now consider $\Omega$ as an open connected set and suppose aditionally that the sequence $\{ f_n \}_{n = 0}^{\infty}$ consists of holomorphic functions. Then for every $z_0 \in \Omega$ there exists an euclidian ball contained in $\Omega$. In view of the Cauchy representation theorem for the holomorphic functions $f_n$
\begin{displaymath}
f_n(\zeta) = \frac{1}{2 \pi \hat{i}} \int\limits_{| z - z_0 |= r} \frac{f_n(z)}{\widetilde{z - \zeta}} d \tilde z,
\end{displaymath}
where every $\zeta$ is such that $|\zeta - z_0| < r$.  In the restriction of every $f_n$ to the closed set $| \zeta - z_0 | \leq \rho < r$ we get
$ r - \rho \leq | z - \zeta | $ which implies 
$(r - \rho)\leq K_2(\alpha, \beta) | z - \zeta |_{(\alpha,\beta)}$.
So we get the following estimation
\begin{displaymath}
\big|\int\limits_{| z - z_0 | = r} \frac{f_n(z) - f(z)}{\widetilde{z - \zeta}} d \tilde z \big|_{(\alpha,\beta)}
\end{displaymath}
\begin{displaymath}
\leq \sup_{| z - z_0 |= r}  |f_n(z) - f(z))|_{(\alpha, \beta)} \frac{1}{  K_1(\alpha, \beta) (r - \rho)} K_2(\alpha, \beta) 2 \pi r.
\end{displaymath}
As the set of restriction is closed the
$f_n$ converges uniformly to $f$ in this compact set
and the right hand-side of the previous expression can be made arbitrarily small. 
We then conclude that for every $\zeta$ such that  $|\zeta - z_0| \leq \rho < r$
\begin{displaymath}
f(\zeta) = \lim_{n \to \infty} f_n(\zeta) =  \lim_{n \to \infty} \frac{1}{2 \pi \hat{i}} \int\limits_{| z - z_0 | = r} \frac{f_n(z) }{\widetilde{z - \zeta}} d \tilde z =  \frac{1}{2 \pi \hat{i}} \int\limits_{| z - z_0 | = r} \frac{f(z) }{\widetilde{z - \zeta}} d \tilde z.
\end{displaymath}
The right hand-side is clearly holomorphic and we conclude that  $f(\zeta)$ is (locally) holomorphic on every closed ball centered in $z_0$. As the point $z_0$ was arbitrarily chose then $f$ is holomorphic in all $\Omega$. We have proved the following:
\begin{theorem}
Suppose $f_n$ is holomorphic in $\Omega$ and the sequence  $\{ f_n \}_{n = 0}^{\infty}$ converges to a limit function $f$ and this convergence is uniform in every subcompact subset of $\Omega$. Then $f$ is holomorphic in $\Omega$.
\end{theorem}

\section{Solution of the initial value problem}

We are now in position to restate the main result of \cite{tut} in the more general context of holomorphy with respect to the structure polynomial $X^2 + \beta X + \alpha$, namely, in order to solve the initial value problem (\ref{unouno}) - (\ref{unocuatro}) we consider an exhaustion of the bounded domain $\Omega$ by a family of subdomains $\Omega_s$, $0 < s < s_0$. The exhaustion can be chosen in such a way that each point $x$ of $\Omega$ lies on the boundary
$\partial \Omega_{s(x)}$ of a uniquely determined domain $\Omega_{s(x)}$ of the exhaustion. Then $s_0 - s(x)$ is a measure of the 
distance of a point $x$ of $\Omega$ from the boundary $\Omega_{s(x)}$. 
Let $\mathcal{B}$ be the Banach space of functions holomorphic in $\Omega_s$ and continuous in $\overline{\Omega_s}$. 
Then the $\mathcal{B}_s$, $0 < s < s_0$ form a scale of Banach spaces. 
Now we rewrite the complex version of the given initial value problem as an abstract operator equation in the scale $\mathcal{B}_s$, 
then we obtain the initial value problem of type
\begin{eqnarray*}
\partial_t w(t, z) & = & \textbf{L} w(t, z) \\
w(0, z) & = & \rho(z),
\end{eqnarray*}
where $w(t, z) = u(t, z) + i v(t, z)$, $z= x + iy$, $t$ means the time, $\textbf{L}$ is defined by (\ref{opeL})
and $\rho(z)= \phi(z) + i \psi(z)$ with $\phi$, $\psi$ given by (\ref{unotres}) and (\ref{unocuatro}) respectively.

In view of the interior estimate of first order for the complex derivative of a holomorphic function and the Weirstrass approximation theorem, a generalized complex abstract Cauchy-Kovaleskaya Theorem is applicable ( see chapter 17 in \cite{treves}). 
Thus, the following theorem has been proved:
\begin{theorem}
The space of holomorphic functions is associated to the system (\ref{unouno}) and (\ref{unodos}) if and only if
the coefficients of the system are given by Lemma 2. For such systems each initial value problem (\ref{unotres}), (\ref{unocuatro})
is solvable, where $\rho = \phi + i \psi $ is an arbitrary holomorphic funtion in $z$ and
the initial functions $\phi$ and $\psi$ satisfy the generalized Cauchy-Riemann system with respect 
to the structure polynomial $X^2 + \beta X + \alpha$ with $4 \alpha - \beta^2 > 0$. 
Moreover the solution writting in the form $w(t, z) = u(t,z) + i v(t, z)$ is holomorphic in z for each t and
exists in the time-interval $0 \leq t \leq h(s_0 - s)$ if $h$ is 
sufficiently small and $(x,y)$ belongs to $\Omega_s$, where the subdomain $\Omega_s$ form an exhaustion of $\Omega$.
\end{theorem}   

\begin{remark}
Since the interior estimate depends on the distance of a subset from the boundary, then
the shorter the distance of the point in $\Omega$
from the boundary of $\Omega$, the shorter the time interval in which the solution exists.
In other words, the solution exists in a conical domain in the $t, x$-space.
\end{remark}

\end{document}